\documentclass[a4paper]{article}

\usepackage{graphicx}
\usepackage{latexsym}
\usepackage{amsmath}
\usepackage{amssymb}

\newenvironment{itemise}{\begin{itemize}}{\end{itemize}}
\newenvironment{centre}{\begin{center}}{\end{center}}

\allowdisplaybreaks[4]

\newcommand{\brak}[1]{\ensuremath{\left[#1\right]}}

\newcommand{\paren}[1]{\ensuremath{\left(#1\right)}}

\newcommand{\goesto}{\rightarrow}

\newcommand{\ints}{\mathbb{Z}}

\newcommand{\ord}[1]{\ensuremath{O\!\paren{#1}}}

\newcommand{\R}{\ensuremath{R_{\alpha}}}
\newcommand{\X}{\ensuremath{X_{\alpha}}}
\newcommand{\I}{\ensuremath{I_{\alpha}}}
\newcommand{\rx}{\ensuremath{\rho_{\alpha}}}
\newcommand{\fl}{\ensuremath{\mu_{\alpha}}}
\newcommand{\cs}{\ensuremath{\sigma_{\alpha}}}
\newcommand{\pc}{\ensuremath{P_{\circ}}}
\newcommand{\rbrch}{\ensuremath{R_{\mathrm{branch}}}}
\newcommand{\xbrch}{\ensuremath{X_{\mathrm{branch}}}}
\newcommand{\pbrch}{\ensuremath{P_{\mathrm{branch}}}}
\newcommand{\zbase}{\ensuremath{Z_{\mathrm{base}}}}
\newcommand{\sbase}{\ensuremath{S_{\mathrm{base}}}}
\newcommand{\vnom}{\ensuremath{V_{\mathrm{nom}}}}

\begin{document}

\title{Power Flows with Flat Voltage Profiles: an Exact Approach}
\author{Anthony B.\ Morton}
\date{July 2022}
\maketitle

\section{Background: the DC Power Flow conditions}
\label{sec:dcflow}

Analysis of alternating--current electrical networks with prescribed power flows is complicated by the fact that network quantities are defined by two degrees of freedom, amplitude and phase displacement.
This greatly complicates the mathematics compared with that of DC circuits and of AC circuits with prescribed currents.

When the focus for analysis is on the flows of active power and on the phase angle relationships between network voltages, it is common to adopt an idealised approach often described as \emph{DC power flow}.
Despite appearances, this does not mean solving a DC network with the same power flows.
Rather, the `DC' is a fossilised reference to DC network analysers, a type of mid-20th century analogue computer that embodied this approach to network solution in discrete electronics.

The `DC power flow' idealisations for an AC power system are as follows:
\begin{itemise}
\item Network branches are assumed lossless or practically lossless, with series reactances $X_{\alpha}$.
\item Voltage magnitudes are assumed equal or nearly equal to nominal values, thus equal to 1 in a per-unit system based on nominal voltage $\vnom$ and an arbitrarily chosen power base $\sbase$.
\end{itemise}
With these idealisations, the full AC power flow equations can be replaced with a simplifed set based on the familiar power-angle formula for a lossless branch
\begin{equation}
P_{\alpha} = \frac{|V_j| |V_k|}{X_{\alpha}} \sin\paren{\delta_j - \delta_k}
\label{eq:pang}
\end{equation}
where $P_{\alpha}$ is the active power flow on branch $\alpha$ (nominally directed from bus $j$ to bus $k$), $X_{\alpha}$ is the branch reactance, $|V_j|$ and $|V_k|$ the terminal bus voltage magnitudes, and $\delta_j$, $\delta_k$ the corresponding voltage phase angles.
Since by assumption $|V_j| = |V_k| = 1$, this further simpifies to
\begin{equation}
P_{\alpha} = B_{\alpha} \sin\paren{\delta_j - \delta_k}
\label{eq:pang1}
\end{equation}
where $B_{\alpha} = 1 / X_{\alpha}$ is the branch susceptance.
Often the approximation is carried one step further, by supposing the phase difference small enough that the sine nonlinearity can be disregarded, thus
\begin{equation}
P_{\alpha} = B_{\alpha} \paren{\delta_j - \delta_k}.
\label{eq:pang2}
\end{equation}
Whether equations (\ref{eq:pang1}) or (\ref{eq:pang2}) are used to represent the network, the state of the idealised system is described by the phase angle $\delta_k$ at each busbar, and the power flow $P_{\alpha}$ on each branch.
In the approximate form (\ref{eq:pang2}), the network equations can be physically modelled by DC circuit branches with conductances $B_{\alpha}$, the DC node voltages serving as analogues of the angles $\delta_k$ and DC branch current (not power) as the analogue of $P_{\alpha}$.
This was the principle behind the DC network analyser.

Neither of the above idealising assumptions is truly valid in practical AC systems, especially that of lossless branches.
All network branches have positive resistance and this can often be significant relative to the branch reactance.
The complication, when nonzero resistances are considered, is that active power losses are determined by both the active and reactive power flows in the branch.

The assumption of uniform nominal or `flat' voltages is also contradicted in practical scenarios, but is arguably more defensible than the zero-resistance assumption, at least for steady-state analysis of transmission networks where relatively abundant voltage controls ensure magnitudes remain within a few percent of nominal values.

This note provides methods for analysing AC power flows when the lossless assumption is discarded but the voltage profile is still assumed `flat' in the sense that all voltages have 1pu magnitude.
By moderating the level of idealisation, these methods can yield results that are more practically relevant in cases where network losses are non-negligible.
They can also assist in demonstrating that theoretical results, obtained for lossless networks, are robust to nonzero resistance values.
One such application is to the existence of circulating flows in wide-area networks \cite{jk:lfiracps}, discussed toward the end of this note.

\section{Definitions}
\label{sec:defs}

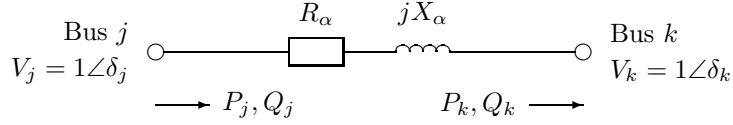
\begin{figure}
\begin{centre}
\begin{picture}(260,40)
\put(50,20){\circle{6}}
\put(40,27){\makebox(0,0)[r]{Bus $j$}}
\put(40,13){\makebox(0,0)[r]{$V_j = 1 \angle \delta_j$}}
\put(53,20){\line(1,0){47}}
\put(50,0){\vector(1,0){20}}
\put(75,0){\makebox(0,0)[l]{$P_j, Q_j$}}
\put(100,15){\framebox(20,10){}}
\put(110,30){\makebox(0,0)[b]{\R}}
\put(120,20){\line(1,0){20}}
\multiput(142.5,20)(5,0){4}{\oval(5,5)[t]}
\put(150,30){\makebox(0,0)[b]{$j\X$}}
\put(160,20){\line(1,0){47}}
\put(210,20){\circle{6}}
\put(220,27){\makebox(0,0)[l]{Bus $k$}}
\put(220,13){\makebox(0,0)[l]{$V_k = 1 \angle \delta_k$}}
\put(190,0){\vector(1,0){20}}
\put(185,0){\makebox(0,0)[r]{$P_k, Q_k$}}
\end{picture}
\end{centre}
\caption{A single branch $\alpha$ in an AC network, with active power flow from bus $j$ to bus $k$.}
\label{fig:branch}
\end{figure}

Figure \ref{fig:branch} depicts a generic branch $\alpha$ within an AC network, having series impedance $\R + j \X$ with $\X > 0$ and $\R \geq 0$.
The relative size of the resistance can be measured by the ratio
\begin{equation}
\rx = \frac{\R}{\X} \geq 0.
\label{eq:rho}
\end{equation}
The terminal buses (nodes) are labelled $j$ and $k$, and for present purposes the voltage magnitude is assumed to be nominal (1pu) at both buses.
The voltage phase angles are denoted $\delta_j$ and $\delta_k$ respectively.

The complex power flow into the branch at bus $j$ (`source end') is defined as $P_j + jQ_j$, with $P$ denoting active power and $Q$ denoting reactive power.
Similarly, the complex power flow out of the branch at bus $k$ (`receiving end') is denoted $P_k + jQ_k$.
For definiteness it will be assumed that the orientation of the branch from $j$ to $k$ matches the direction of active power flow, and so
\begin{equation}
P_j \geq P_k \geq 0
\label{eq:psense}
\end{equation}
with $P_j = P_k$ when the branch is lossless ($\R = 0$).

The current flow in the branch will be denoted \I, with magnitude $|\I|$, in terms of which the active and reactive power consumption satisfy the identities
\begin{equation}
P_j - P_k = \R |\I|^2 \geq 0, \qquad Q_j - Q_k = \X |\I|^2 \geq 0,
\label{eq:pqsunk}
\end{equation}
from which also follows the identity
\begin{equation}
\fl = \X P_j - \R Q_j = \X P_k - \R Q_k.
\label{eq:mu}
\end{equation}
The \emph{flow coefficient} $\fl$ will play an important role in the power flow solution below.

Lastly, the notation $\cs$ will denote the (dimensionless) ratio between the reactive power consumption on the branch and the receiving end active power:
\begin{equation}
\cs = \frac{\X |\I|^2}{P_k} \qquad \text{(when $P_k > 0$, else $\cs = 0$).}
\label{eq:sigma}
\end{equation}
This is denoted the branch \emph{coefficient of support} and will also prove useful below.

\section{Exact solution for a branch}
\label{sec:solution}

From Ohm's law together with the fundamental complex power and voltage identities
\begin{equation}
V_j \I^* = P_j + jQ_j, \qquad V_k \I^* = P_k + jQ_k, \qquad
V_j V_k^* = |V_j| |V_k| \exp j\!\paren{\delta_j - \delta_k},
\label{eq:vcplx}
\end{equation}
can be derived two key formulae that provide the exact power flow solution for the branch $\alpha$ of Figure \ref{fig:branch}.
The first is the identity relating the flow coefficient to the phase angle shift
\begin{equation}
\fl = \X P_j - \R Q_j = \X P_k - \R Q_k = |V_j| |V_k| \sin\paren{\delta_j - \delta_k}
\label{eq:flow}
\end{equation}
which generalises formula (\ref{eq:pang}) to branches with arbitrary series resistances.
The second is a quadratic equation relating the (squared) magnitudes of the endpoint voltages, written here for the receiving end $k$:
\begin{equation}
|V_k|^4 - \brak{|V_j|^2 - 2 \paren{\R P_k + \X Q_k}} |V_k|^2
   + \paren{\R^2 + \X^2} \paren{P_k^2 + Q_k^2} = 0.
\label{eq:vmag}
\end{equation}
Equation (\ref{eq:vmag}) is a useful tool for determining $|V_k|$ in situations where $|V_j|$, $P_k$ and $Q_k$ are known.
Together with (\ref{eq:flow}) it provides the general branch flow solution when any three of the five quantities $|V_j|$, $|V_k|$, $P_k$, $Q_k$ and $\delta_j - \delta_k$ are known.

Under the `flat voltage' assumption $|V_j| = |V_k| = 1$, formula (\ref{eq:vmag}) leads directly to the condition
\begin{equation}
\paren{\R^2 + \X^2} \paren{P_k^2 + Q_k^2} = -2 \paren{\R P_k + \X Q_k}.
\label{eq:flatpq}
\end{equation}
For a given $P_k$, this yields a quadratic equation for $Q_k$
\begin{equation}
\paren{\R^2 + \X^2} Q_k^2 + 2 \X Q_k + \brak{2 \R P_k + \paren{\R^2 + \X^2} P_k^2} = 0
\label{eq:flatpqquad}
\end{equation}
with solutions
\begin{equation}
Q_k = \frac{1}{\R^2 + \X^2}
   \brak{-\X \pm \sqrt{\X^2 - 2 \R \paren{\R^2 + \X^2} P_k - \paren{\R^2 + \X^2}^2 P_k^2}}.
\label{eq:flatq1}
\end{equation}
Observe that with $\X > 0$ and $\R, P_k \geq 0$ by assumption, both solutions for $Q_k$ are negative or zero (the latter only when $P_k = 0$).
The solution with the positive sign is closer to zero, while the other is of order $1 / \X$ (and lies on an `inverted branch' of the network $Q$--$V$ curve).
Taking the former, therefore, as the `practical' solution for $Q_k$ and invoking the ratio (\ref{eq:rho}) gives
\begin{equation}
Q_k = - \frac{1 / \X}{1 + \rx^2}
   \brak{1 - \sqrt{1 - 2 \rx \paren{1 + \rx^2} \paren{\X P_k} - \paren{1 + \rx^2}^2 \paren{\X P_k}^2}}.
\label{eq:flatq}
\end{equation}
Except in the quiescent condition $P_k = 0$, the value $Q_k$ given by (\ref{eq:flatq}) is strictly negative.
That is to say, maintaining a flat voltage profile necessitates that there be an injection of reactive power $Q_k$ at the receiving end of the branch, in the contrary direction to the active power $P_k$.

Aside from confirming $Q_k$ is negative, (\ref{eq:flatq}) does not give immediate insight into the quantity required, although it points to a breakdown in the solution (voltage collapse) when $P_k$ is of order $1/\X$.
In most practical situations the magnitude of $P_k$ is a good deal smaller, and it is helpful to approximate the square root using the Taylor expansion
\begin{equation}
\sqrt{1 - x} = 1 - \frac{x}{2} - \frac{x^2}{8} - \frac{x^3}{16} - \ord{x^4}
\label{eq:sqrtapprox}
\end{equation}
which leads to
\begin{multline}
\sqrt{1 - 2 \rx \paren{1 + \rx^2} \paren{\X P_k} - \paren{1 + \rx^2}^2 \paren{\X P_k}^2} \\
   = 1 - \rx \paren{1 + \rx^2} \paren{\X P_k} - \frac{\paren{1 + \rx^2}^3}{2} \paren{\X P_k}^2
      - \frac{\rx \paren{1 + \rx^2}^4}{2} \paren{\X P_k}^3 - \ord{\paren{\X P_k}^4}.
\label{eq:flatqq}
\end{multline}
Substituting in (\ref{eq:flatq}) and simplifying yields a more practical formula for $Q_k$ when $\X P_k \ll 1$:
\begin{equation}
Q_k = - P_k \brak{\rx + \frac{\paren{1 + \rx^2}^2}{2} \paren{\X P_k}
   + \frac{\rx \paren{1 + \rx^2}^3}{2} \paren{\X P_k}^2 + \ord{\paren{\X P_k}^3}}.
\label{eq:flatqp}
\end{equation}
This formula reveals a qualitative difference between the reactive support requirements in a near-lossless branch relative to one with substantial resistance.
In the lossless or near-lossless case, $\rx$ is practically zero and so the leading term vanishes; the required reactive power is $Q_k \approx - (\X / 2) P_k^2$ to leading order.
On the other hand, when $\rx$ is larger the leading term dominates, so that $Q_k \approx -\rx P_k$ to leading order, and scales linearly with active power flow rather than quadratically.

Having solved for $Q_k$ it is straightforward to determine the current $\I$, sending-end reactive power $Q_j$, and the phase displacement $\delta_j - \delta_k$.
To obtain practically useful formulae for these quantities, however, it is useful to introduce a dimensionless parameter that can temporarily set aside purely mathematical complications.

\section{The Coefficient of Support}
\label{sec:coeff}

The \emph{coefficient of support} $\cs$ is defined by (\ref{eq:sigma}) above:
\[\cs = \frac{\X |\I|^2}{P_k} \qquad \text{(when $P_k > 0$, else zero)}.\]
To recap, if $P_k$ is the power flow at the `receiving' end of the branch $\alpha$, then by definition $\Delta Q = \cs P_k$ is the total reactive power consumed in the branch reactance.
Reactive power is a conserved quantity, so there must be a compensating net injection of reactive power at the terminal buses to `support' the power flow $P_k$.

Assuming $\X > 0$, the reactive power requirement is always positive when $P_k > 0$.
If $P_k = 0$, then there is clearly no current flow on the branch and so the reactive power consumption is also zero: it will be seen below that in fact $\cs \goesto 0$ as $P_k \goesto 0$, which motivates extending the definition (\ref{eq:sigma}) to set $\cs = 0$ when $P_k = 0$.

The coefficient $\cs$ can be expressed explicitly in terms of $P_k$ and the branch impedance alone, using the solution in the previous section.
The identity $V_k \I^* = P_k + jQ_k$ immediately gives
\begin{equation}
|\I|^2 = \frac{1}{|V_k|^2} \paren{P_k^2 + Q_k^2} = P_k^2 + Q_k^2
\label{eq:ibranch}
\end{equation}
and combining with the condition (\ref{eq:flatpq}) for a flat voltage profile gives
\begin{equation}
|\I|^2 = P_k^2 + Q_k^2 = -\frac{2}{\R^2 + \X^2} \paren{\R P_k + \X Q_k}.
\label{eq:ibranch2}
\end{equation}
(Note this implies that $Q_k$ is both negative and large enough that $\R P_k + \X Q_k$ is also negative.)
It follows that
\begin{equation}
\cs = - 2 \frac{\R \X}{\R^2 + \X^2} - 2 \frac{\X^2}{\R^2 + \X^2} \frac{Q_k}{P_k}
   = \frac{2}{1 + \rx^2} \paren{\frac{(-Q_k)}{P_k} - \rx},
\label{eq:cs}
\end{equation}
where $Q_k$ is given in terms of $P_k$ and the branch impedance by equation (\ref{eq:flatq}).
The exact formula for $\cs$ that results is cumbersome, but when the asymptotic formula (\ref{eq:flatqp}) for $Q_k$ is used instead, a much nicer expression results:
\begin{equation}
\cs = \paren{1 + \rx^2} \paren{\X P_k} + \rx \paren{1 + \rx^2}^2 \paren{\X P_k}^2
   + \ord{\paren{\X P_k}^3}.
\label{eq:csp}
\end{equation}
It is evident that this expression for $\cs$ can in principle be continued as a power series in the quantity $\X P_k$, the coefficients being polynomials in $\rx$.
To leading order $\cs$ is just $(\X + \R^2 / \X) P_k$, and if the branch is lossless than $\cs$ is just $\X P_k$ with residual error of order $(\X P_k)^3$.

The power series (\ref{eq:csp}) also hints at $\cs$ being a strictly increasing function of $P_k$, as intuition might suggest.
This property can be confirmed by direct calculation using (\ref{eq:cs}) and (\ref{eq:flatq}), which reveals that
\begin{equation}
\frac{\partial \cs}{\partial P_k} = \frac{\cs}{P_k \sqrt{\Delta}} > 0 \qquad \paren{P_k > 0}
\label{eq:csderiv}
\end{equation}
where $\Delta$ is the discriminant under the square root in (\ref{eq:flatq}).
When $P_k \goesto 0$ one may use L'H\^{o}pital's rule to find
\begin{equation}
\lim_{P_k \goesto 0} \frac{\partial \cs}{\partial P_k}
   = \frac{2}{1 + \rx^2} \lim_{P_k \goesto 0} \paren{\frac{1}{2} \frac{\partial^2 (-Q_k)}{\partial P_k^2}}
   = \paren{1 + \rx^2} \X > 0,
\label{eq:csderiv0}
\end{equation}
agreeing with the power series (\ref{eq:csp}) and showing that $\cs$ is strictly increasing with $P_k \geq 0$.

With a little more effort it may also be shown that provided $P_k > 0$, $\cs$ is strictly increasing with the ratio $\rx$, all else being equal (in other words, strictly increasing with $\R$ for fixed $\X$).
Starting again from (\ref{eq:flatq}) and some algebraic manipulations, one may show that
\begin{equation}
\frac{\partial (-Q_k)}{\partial \rx} = \frac{2 \rx}{\paren{1 + \rx^2} \sqrt{\Delta}} \paren{-Q_k}
      + \frac{1 - \rx^2}{1 + \rx^2} \frac{P_k}{\sqrt{\Delta}}
   = \frac{P}{\sqrt{\Delta}} \paren{1 + \rx \cs},
\label{eq:dqdr}
\end{equation}
from which it can be shown using (\ref{eq:cs}) that
\begin{equation}
\frac{\partial \cs}{\partial \rx} = \frac{2 \X}{\sqrt{\Delta}} \paren{1 + \rx \cs} \paren{-Q_k} > 0 \qquad
\text{when } P_k > 0.
\label{eq:csderivr}
\end{equation}
Of course when $P_k = 0$, then $\cs = 0$ has no dependence on $\rx$.

It should likewise come as little surprise that $\cs$ is strictly increasing with $\X$, when all other parameters are held constant (including the ratio $\rx$).
Direct calculation using (\ref{eq:cs}) and (\ref{eq:flatq}) confirms that indeed
\begin{equation}
\frac{\partial \cs}{\partial \X} = \frac{\cs}{\X \sqrt{\Delta}} > 0 \qquad \text{when } P_k > 0.
\label{eq:csderivx}
\end{equation}
As a consequence of all the above, it may be seen that if two of the three parameters $P_k$, $\X$ and $\rx$ are fixed, there is a one-to-one mapping between $\cs$ and the third parameter in the set, such that the two variables increase or decrease together.

The coefficient of support yields simple exact formulae for all quantities associated with the branch power flow solution, given $P_k$ and the branch impedance, as follows.
\begin{alignat}{2}
\text{Branch current:} &\quad& |\I| &= \sqrt{\frac{\cs P_k}{\X}}
\label{eq:csi} \\
\text{Power losses:} &\quad& P_j - P_k = \R |\I|^2 &= \rx \cs P_k
\label{eq:csloss} \\
\text{Receiving end Q:} &\quad& Q_k &= \paren{- \frac{1 + \rx^2}{2} \cs - \rx} P_k
\label{eq:csqk} \\
\text{Sending end Q:} &\quad& Q_j = Q_k + \cs P_k
   &= \paren{\phantom{-} \frac{1 - \rx^2}{2} \cs - \rx} P_k
\label{eq:csqj} \\
\text{Sending end P:} &\quad& P_j &= \paren{1 + \rx \cs} P_k
\label{eq:cspj} \\
\text{Flow coefficient:} &\quad& \fl = \X P_k - \R Q_k
   &= \X P_k \paren{1 + \rx^2} \paren{1 + \frac{\rx \cs}{2}}
\label{eq:csfl} \\
\text{Phase shift:} &\quad& \delta_j - \delta_k &= \arcsin \fl.
\label{eq:csdelta}
\end{alignat}
Formulae (\ref{eq:csqk}) and (\ref{eq:csqj}), specifically, indicate how reactive power support must be configured at the two ends of the branch in order to maintain the active power flow $P_k$ under a flat voltage profile.
For a lossless branch, this amounts to a symmetric injection $Q = (\cs / 2) P$ at both ends.
With nonzero resistance, there must be an additional $Q$ flow superimposed on this symmetric injection, contrary to the direction of $P$ flow and with magnitude $(\rx + \rx^2 \cs / 2) P_k$.

Based on the fact observed above that $\cs$ is strictly increasing with $P_k \geq 0$, one may also confirm that the current $|\I|$, the power losses, the flow coefficient and the phase shift are also all strictly increasing functions of $P_k \geq 0$.

There is of course a limit to the power flow $P_k$ that can be sustained in this manner: as $P_k$ increases, the discriminant $\Delta$ under the square root in (\ref{eq:flatq}) decreases until eventually it reaches zero, beyond which point there is no solution for $Q_k$.
Equating $\Delta$ to zero reveals that the limiting values for $P_k$ and $Q_k$ are
\begin{equation}
P_k = \frac{1}{\X} \cdot \frac{\sqrt{1 + \rx^2} - \rx}{1 + \rx^2}, \qquad \text{and} \qquad
Q_k = - \frac{1}{\X} \cdot \frac{1}{1 + \rx^2}.
\label{eq:limpq}
\end{equation}
Formulae (\ref{eq:cs}) and (\ref{eq:csfl}) and a little algebra than show that under the same conditions
\begin{equation}
\cs = \frac{2}{\sqrt{1 + \rx^2}} \qquad \text{and} \qquad \fl = \frac{1}{\sqrt{1 + \rx^2}}.
\label{eq:limcs}
\end{equation}
These identities have notable practical implications for the limiting power flow condition: first, that at the limit the reactive power consumption on the branch can be up to twice the active power flow; and second, that the phase shift in bus voltages at the limit is equal to the `impedance angle'
\begin{equation}
\phi = \arcsin\frac{\X}{\sqrt{\R^2 + \X^2}}.
\label{eq:zang}
\end{equation}
In particular, if $\R > 0$ then the maximum phase shift is strictly less than $\pi/2$.
When $\R = 0$ on the other hand one has $P_k = -Q_k = 1 / \X$, $\cs = 2$, $\fl = 1$ and the phase shift is $\pi/2$ at the limit.

\section{The power--angle relation with flat voltages}
\label{sec:pang}

As noted at the outset, the `DC power flow' approach to network solution hinges on the relation (\ref{eq:pang}) between phase displacement and active power flow in a lossless branch:
\[\sin\paren{\delta_j - \delta_k} = \X P_{\alpha}. \qquad \paren{|V_j| = |V_k| = 1, \R = 0}\]
It was seen above that after factoring in resistance, this must be modified to equation (\ref{eq:csfl}):
\[\sin\paren{\delta_j - \delta_k} = \fl = \X P_k \paren{1 + \rx^2} \paren{1 + \frac{\rx \cs}{2}},\]
where $P_k \leq P_j$ is the receiving end power, $\rx$ is the $R$-to-$X$ ratio, and the coefficient of support $\cs$ increases monotonically with any one of $(P_k, \X, \rx)$ when the others are held fixed.

The monotonicity properties of $\cs$ assist greatly in understanding the \emph{qualitative} behaviour of the relation (\ref{eq:csfl}) between power flow and angle displacement.
To begin with, it asserts that all else being equal, an increase (reduction) in any of the parameters $P_k$, $\X$ or $\R$ will also increase (reduce) $\fl$ and, consequently, the angle displacement.
Importantly, it also asserts that if the value of $\fl$ is fixed, then an increase (reduction) in one of these same parameters requires a corresponding reduction (increase) in one of the others to keep $\fl$ constant.

Equation (\ref{eq:csfl}) is however of limited assistance in obtaining \emph{quantitative} solutions for $\fl$ or (especially) $P_k$.
This requires working with an explicit formula for $\cs$ in terms of $(P_k, \X, \rx)$, which can be obtained using (\ref{eq:cs}) from the previous section and substituting (\ref{eq:flatq}) for $Q_k$.
When the resulting formula for $\cs$ is used in (\ref{eq:csfl}), the result is
\begin{equation}
\fl = \X P_k + \frac{\rx}{1 + \rx^2} \paren{1
   - \sqrt{1 - 2 \rx \paren{1 + \rx^2} \paren{\X P_k} - \paren{1 + \rx^2}^2 \paren{\X P_k}^2}}.
\label{eq:flowxp}
\end{equation}
This formula remains cumbersome, but simplifies substantially when rearranged to clear the square root.
A subtle point emerges when doing so however: observe that the nonnegativity of the square root in (\ref{eq:flowxp}) implies that $\fl$ and $P_k$ must satisfy the inequality
\begin{equation}
\fl - \X P_k \leq \frac{\rx}{1 + \rx^2}.
\label{eq:flowpcheck}
\end{equation}
(This inequality ultimately stems from having excluded at the outset the `impractical' solution for the reactive power $Q_k$.)
Since (\ref{eq:csfl}) already implies that $\fl \geq \X P_k$, one now has two-sided bounds for the difference between $\fl$ and $\X P_k$ that arises when $\rx > 0$.
Furthermore, given that the quantity under the square root in (\ref{eq:flowxp}) is the same discriminant $\Delta$ from (\ref{eq:flatq}), one can readily characterise the conditions leading to these bounds: the difference $\fl - \X P_k$ is a minimum (zero) for a lossless branch, and is a maximum (\ref{eq:flowpcheck}) in the limiting flow condition given by (\ref{eq:limpq}) and (\ref{eq:limcs}).

Isolating the square root in (\ref{eq:flowxp}), squaring and simplifying now yields a quadratic equation linking $P_k$ and $\fl$:
\begin{equation}
\paren{1 + \rx^2} \fl^2 - 2 \rx \fl + 2 \paren{\rx - \fl} \paren{1 + \rx^2} \paren{\X P_k}
   + \paren{1 + \rx^2}^2 \paren{\X P_k}^2 = 0.
\label{eq:flp}
\end{equation}
When this is rewritten to complete the square in $P_k$ as
\begin{equation}
\brak{\paren{1 + \rx^2} \paren{\X P_k} + \paren{\rx - \fl}}^2 + \paren{\fl^2 - 1} \rx^2 = 0,
\label{eq:flpp}
\end{equation}
it can be solved directly for $P_k$ given $\fl$.
This would ordinarily yield two solutions for $P_k$, but one may verify that only one of these actually satisfies the inequality (\ref{eq:flowpcheck})---the other solution enters when squaring both sides to clear the square root in (\ref{eq:flowxp}), and does not actually satisfy (\ref{eq:flowxp}) itself.
Accordingly, one has the unique solution for $P_k$ given $\fl$
\begin{equation}
P_k = \frac{1 / \X}{1 + \rx^2} \brak{\fl - \rx \paren{1 - \sqrt{1 - \fl^2}}}.
\label{eq:pfl}
\end{equation}
Note first that for $\rx = 0$ this reproduces the lossless solution (\ref{eq:pang}), in the form $\X P_k = \fl$.
One may also verify that when $\fl$ takes its limiting value $1 / \sqrt{1 + \rx^2}$, formula (\ref{eq:pfl}) reduces to formula (\ref{eq:limpq}) for $P_k$ as it should.
Finally, it may be verified that
\begin{equation}
\frac{\partial P_k}{\partial \fl} = \frac{1 / \X}{1 + \rx^2} \paren{1 - \frac{\rx \fl}{\sqrt{1 - \fl^2}}},
\label{eq:pflderiv}
\end{equation}
and this is strictly positive for $\fl < 1 / \sqrt{1 + \rx^2}$, decreasing to zero when $\fl$ reaches this limiting value.
This confirms that $P_k$ and $\fl$ always increase together when the branch impedance is held constant, hence that the $P_k$ value (\ref{eq:limpq}) at the flow limit is also the largest active power transfer that can be achieved on a branch of fixed impedance, subject to the flat voltage profile.

\section{Circulating flows and winding numbers in large networks}
\label{sec:circ}

The preceding analysis for a single network branch can be applied to each of the branches making up a wider network, with varying power flows $P_k$ on each branch and injections of active and reactive power at the network buses.
Equation (\ref{eq:csfl}) and its converse (\ref{eq:pfl}) provide the necessary linkage between each branch power flow, the impedance, and the displacement between the terminal bus voltage phase angles.
In particular, they show that when branch resistance is taken into account, the angle displacement is always larger than suggested by the formula $\sin(\delta_j - \delta_k) = \X P_k$ for lossless branches.
The difference between $\fl = \sin(\delta_j - \delta_k)$ and $\X P_k$ ranges from zero (when the branch is lossless) up to the limit given by equation (\ref{eq:flowpcheck}) in terms of the $R/X$ ratio $\rx$.

The most immediate application is to a network with `string' topology, where successive branches are connected end-to-end, and the sending power $P_j$ on one branch is the sum of receiving power $P_k$ for the previous branch and any additional power injected at the intermediate bus.
Assuming each $P_k > 0$, the voltage phase angles $\delta_k$ at successive busbars form a strictly decreasing sequence, the lag being determined at each stage by the power flow $P_k$ and branch impedance.
Reactive power injection is also required at each intermediate busbar in order to meet the support requirement $\Delta Q = \cs P_k > 0$ on each branch and maintain the flat voltage profile.

It turns out to be possible to `close the loop' on such a network, identifying an upstream bus with a remote downstream bus to form a ring network with a positive circulating power flow.
Let such a ring network comprise $n$ buses (numbered $0$ through $n - 1$) and $n$ branches (numbered $1$ through $n$), with branch $k$ linking buses $k - 1$ and $k$, and bus indices $n$ and $0$ identified.
For $1 \leq k \leq n$ it is assumed that $P_k$, the receiving-end power on the $k$th branch, is strictly positive.
Then one may identify the \emph{circulating power} $\pc$ as the minimum of $P_k$ over all indices $k$, since the measured active power flow is uniform in direction and no less than $\pc$ at any point in the ring.

Of course, since $\delta_{k-1} - \delta_k > 0$ for all $k$ when $\pc > 0$ it is not possible for $\delta_0$ and $\delta_n$ to be identical as real numbers, but this is not necessary: as they represent angle variables it is only necessary that they differ by an integer multiple of $2\pi$ \cite{jk:lfiracps}.
That is to say,
\begin{equation}
\sum_{k=1}^n \paren{\delta_{k-1} - \delta_k} = \sum_{\alpha=1}^n \arcsin\fl = 2 m \pi,
\qquad m \in \ints_+.
\label{eq:winding}
\end{equation}
The integer $m$ is termed the \emph{winding number} in the literature, and can be defined for any cycle of branches in networks with node states characterised by angles, on the understanding that the differences $(\delta_{k-1} - \delta_k)$ in (\ref{eq:winding}) are always restricted to the interval $(-\pi,\pi]$ \cite{jhsb:fenontgac}.
As seen above, a nonzero winding number is necessary for circulating power to exist, and so a winding number of zero for all cycles in the network certifies the absence of circulating power flow.

Circulating power flows can in principle exist in any ring network of sufficient size, with (in general) non-uniform impedances $\R + j\X$ and nonuniform power flows, subject only to there being net power injections at each busbar sufficient to maintain the reactive power consumption $\cs P_k$ and active power losses $\rx \cs P_k$ on every branch, for some $P_k \geq \pc$.
However, the `homogenous' case with uniform impedances and power flows is worth considering in further detail as a pointer to the behaviour in the general case.

Suppose then that in a ring network of $n$ branches, every branch has the same impedance $R + jX$ and the same power flow $P = \pc$ at the receiving end, directed uniformly.
(Subscripts are now dropped given all branches are identical.)
To relate this to real-world applications it is useful to recall that these quantities are expressed in a per-unit system, so that if $\rbrch$, $\xbrch$ and $\pbrch$ are the equivalent branch quantities in SI units and $\vnom$ is the nominal voltage at each bus, then
\begin{equation}
R = \frac{\rbrch}{\zbase}, \quad X = \frac{\xbrch}{\zbase}, \quad P = \frac{\pbrch}{\sbase} \qquad
\text{where} \qquad \zbase = \frac{\vnom^2}{\sbase},
\label{eq:perunit}
\end{equation}
with $\sbase$ chosen arbitrarily as a convenient base scale for power flows (often taken as 100MVA in large transmission systems).

Since each branch has the same impedance and power flow, each also has the same coefficient of support $\sigma$ and the same flow coefficient $\mu$.
The winding number condition (\ref{eq:winding}) therefore becomes
\begin{equation}
\mu = \sin\paren{\frac{2m}{n} \pi}, \qquad m \in \ints_+, \qquad
   1 \leq m \leq \left\lfloor \frac{n}{4} \right\rfloor,
\label{eq:ringflow}
\end{equation}
where the restriction on $m$ is to ensure $\delta_j - \delta_k \leq \pi/2$ for each individual branch, and implies in particular that there are at least 4 branches making up the ring.
(In the context of steady state power flows, angle differences greater than $\pi/2$ across a single branch are not practical, as they invert the increasing relationship between active power and angle difference that underpins system stability, and also, much like the second solution to (\ref{eq:flatq1}), lead to impractical reactive power flows.)

Formula (\ref{eq:pfl}) then implies a unique circulating power
\begin{equation}
\pc = \frac{1 / X}{1 + \rho^2}
   \brak{\sin\paren{\frac{2m}{n} \pi} - \rho \paren{1 - \cos\paren{\frac{2m}{n} \pi}}},
\label{eq:ringpc}
\end{equation}
for each admissible winding number $m$, given the impedance of a single branch and the number $n$ of branches.
This is clearest in the lossless case $R = 0$, where the solutions are
\begin{equation}
\pc = \frac{1}{X} \sin\paren{\frac{2m}{n} \pi} \qquad
\text{for} \qquad 1 \leq m \leq \left\lfloor \frac{n}{4} \right\rfloor.
\label{eq:ringlossless}
\end{equation}
The practice is adopted here and below (as in both (\ref{eq:limpq}) and (\ref{eq:pfl}) above) of expressing power quantities as a multiple of $1/X$, which provides a convenient scale for circulating flows.
Note from above that using a scale of $1/X$ for per-unit power is equivalent to using a scale of $\vnom^2 / \xbrch$ for power in engineering units (watts or VA).

For a more quantitative analysis, consider first the limiting values of $\rho = R / X$ and $\pc$ for fixed values of $n$ and $X$.
Recall that from (\ref{eq:limcs}) the maximum feasible value of $\mu$ is $1 / \sqrt{1 + \rho^2}$.
Combined with taking $m = 1$ in (\ref{eq:ringflow}) to get the \emph{smallest} feasible value of $\mu$ for fixed $n$, this prescribes an \emph{upper} bound on $\rho$ to achieve any circulating flow in a ring of $n$ branches.
This upper bound $\rho_{\max}$ is tabulated in Table \ref{tab:rlimit} along with the $\pc$ value at $\rho_{\max}$ from either (\ref{eq:limpq}) or (\ref{eq:ringpc}), the reactive power consumption $\Delta Q = \sigma \pc$ (where $\sigma = 2 \mu$ by (\ref{eq:limcs})), and the active power losses on each branch.
\begin{table}
\begin{centre}
\begin{tabular}{c|c|c|c|c}
Ring & Largest & Limiting flow & $Q$ consumption & $P$ losses \\
size & feasible & $\pc(\rho_{\max})$ & per branch $\sigma \pc$ & per branch $\rho \sigma \pc$ \\
$n$ & $R / X$ ratio & \multicolumn{3}{c}{\small(all in units of $\vnom^2 / \xbrch$)} \\ \hline
4 & 0 & 1 & 2 & 0 \\
5 & 0.3249 & 0.6572 & 1.25 & 0.4061 \\
6 & 0.5774 & 0.4330 & 0.75 & 0.4330 \\
7 & 0.7975 & 0.2944 & 0.4603 & 0.3671 \\
8 & 1 & 0.2071 & 0.2929 & 0.2929 \\
9 & 1.1918 & 0.1504 & 0.1933 & 0.2304 \\
10 & 1.3764 & 0.1123 & 0.1320 & 0.1817
\end{tabular}
\end{centre}
\caption{Limits on $R / X$ and circulating power for ring networks with small numbers of branches.}
\label{tab:rlimit}
\end{table}

It should be stressed the values in Table \ref{tab:rlimit} only pertain to the (normalised) limiting flow at this largest feasible $R/X$ ratio $\rho_{\max}$.
So it should not be too surprising that for this extreme limiting flow, the $Q$ consumption on each branch is often larger than the actual active power flow $\pc$---indeed by (\ref{eq:limcs}) the coefficient of support $\sigma$ at the limit exceeds 1 whenever $\rho < \sqrt{3}$.
Likewise, for networks of 7 branches or more having the largest possible resistance in accordance with $\rho_{\max}$, the $P$ losses on each branch actually exceed $\pc$.
(The threshold case is that of 6 branches, where $\mu = \sqrt{3} / 2$, $\rho_{\max} = 1 / \sqrt{3}$ and $\sigma = \sqrt{3}$, so $\rho \sigma = 1$ and the losses are the same as $\pc$, both $\sqrt{3} / 4$ in per-unit.)

If the actual $R/X$ ratio is smaller than $\rho_{\max}$, flows with larger circulating power $\pc$ are possible.
So for example when $n = 7$, a flow with winding number $m = 1$ must have $\mu = \sin(2\pi/7) = 0.7818$ irrespective of the ratio $\rho = R/X$, but any value $\rho \leq 0.7975$ yields feasible solutions $\pc$ given by equation (\ref{eq:pfl}), right down to $\rho = 0$ where the solution is $\pc = \mu / X = 0.7818 / X$.
It is readily seen from (\ref{eq:pfl}) that $X \pc$ increases monotonically with decreasing $\rho$ when $\mu$ is fixed, from $0.2944$ at $\rho = \rho_{\max}$ to $0.7818$ at $\rho = 0$.
This also follows from equation (\ref{eq:csfl}) due to the monotonicity properties of $\sigma$ found earlier.

It follows that for a fixed reactance $X$ per branch and a fixed number of branches $n$, the values $\pc(\rho_{\max})$ in Table \ref{tab:rlimit} are in fact the \emph{minimum} feasible circulating flows, corresponding to the \emph{maximum} feasible network resistances.
Conversely, the \emph{maximum} feasible circulating flows with a fixed $X$ and $n$ (and winding number 1) are the values given by (\ref{eq:ringlossless}) with $m = 1$ and $R = \rho = 0$.
For values intermediate between these extremes, there will be a ratio $\rho = R/X \in (0, \rho_{\max})$ that achieves this value of circulating power with a winding number of 1.

Larger winding numbers are achievable when $n$ is sufficiently large, as suggested by (\ref{eq:ringflow}).
Each winding number generates a fixed value of $\mu$ by (\ref{eq:ringflow}), which when equated to $1 / \sqrt{1 + \rho^2}$ yields an equivalent $\rho_{\max}$ value.
A continuum of feasible circulating flows $\pc$ again results, corresponding to values of $\rho$ between zero and $\rho_{\max}$.
The endpoints of the feasible continuum for $\pc$ are given by (\ref{eq:ringlossless}) for $\rho = 0$ and by (\ref{eq:limpq}) or (\ref{eq:ringpc}) for $\rho = \rho_{\max}$.

One may further note, as already found above, that whenever $R > 0$ the circulating power $\pc$ is always accompanied by a circulating $Q$ flow in the opposite direction, of magnitude $(\rho + \rho^2 \sigma / 2) \pc$.
This circulating $Q$ flow is additional to the consumption $\Delta Q = \sigma \pc$ on each branch, which must be supplied as a $Q$ injection at each intervening busbar.

\section{Conclusion}
\label{sec:conclu}

Working through the exact AC power flow solution with a `flat' voltage profile provides important insights into the relationship between active and reactive power flows, reactive compensation requirements and phase angle differences when the resistance of network branches is taken into account.
Through use of the `coefficient of support' associated with the power flow on each network branch, the inevitable mathematical complications can be de-emphasised without sacrificing precision.

The simple relation $\sin(\delta_j - \delta_k) = X P$ for a lossless network branch is generalised by recognising that with a flat voltage profile, the correct value on the right hand side of this relation is the flow coefficient $\mu = X P - R Q$.
Solving for $Q$ from first principles, it is found that a positive injection of reactive power is always required to sustain a flat voltage profile with $|V| = 1$pu at both ends of a network branch with positive reactance $X$.
The required injection is readily quantifiable in terms of the coefficient of support $\sigma = X |I|^2 / P$.
For branches with positive resistance one always has $\sin(\delta_j - \delta_k) > X P$, and bounds for the difference $\sin(\delta_j - \delta_k) - X P$ are readily derived.

The generalised analysis has a ready application to the study of circulating power flows that can arise in large ring systems.
For the case of a uniform ring with a circulating flow, the analysis provides bounds on the maximum allowable network resistance and the minimum attainable power flow in terms of the other configuration parameters.
It also reveals that when resistance is taken into account, a real circulating power flow is always accompanied by a contrary flow of reactive power.

\bibliographystyle{plain}
\bibliography{systems}

\end{document}